\documentclass[a4paper]{elsarticle}

\usepackage{amsmath,amsfonts,amssymb,url,theorem}

\usepackage[utf8]{inputenc}

{\theorembodyfont{\slshape}
\newtheorem{proposition}{Proposition}[section]
\newtheorem{lemma}[proposition]{Lemma}

\newtheorem{theorem}[proposition]{Theorem}}

{\theorembodyfont{\upshape}
\newtheorem{remark}[proposition]{Remark}
}

\newcommand\C{{\mathbb C}}
\renewcommand\H{{\mathbb H}}

\newcommand\Q{{\mathbb Q}}

\newcommand\Z{{\mathbb Z}}

\newcommand{\OO}{{\mathcal O}}

\newcommand\calA{{\mathcal A}}

\newcommand\Cl{{\mathrm{Cl}}}
\newcommand\End{{\mathrm{End}}}

\renewcommand\Im{{\mathrm{Im}}}
\renewcommand\Re{{\mathrm{Re}}}

\newcommand\gerD{{\mathfrak D}}

\title{Rational Products of Singular Moduli}

\author[bdx]{Yuri Bilu\fnref{hamot,algant}} 
\author[wits]{Florian Luca\fnref{algant}}
\author[valpo]{Amalia Pizarro-Madariaga\fnref{algant}}

\address[bdx]{Institut de Mathématiques de Bordeaux,
Université de Bordeaux and CNRS,
Talence, France}

\address[wits]{School of Mathematics, University of the Witwatersrand,  South Africa}

\address[valpo]{Departamento de Matemáticas, Universidad de Valpara\'iso, Chile}

\fntext[hamot]{Supported by  the \textsl{Agence National de la Recherche} project ``Hamot'' (ANR 2010 BLAN-0115-01)}

\fntext[algant]{Supported by the ALGANT Scholarship program}

\date\today

\begin{document}

\hfuzz 3pt

\begin{abstract}
We show that with ``obvious'' exceptions the product of two singular moduli cannot be a non-zero rational number. This gives a totally explicit version of André's 1998 theorem on special points for the hyperbolas ${x_1x_2=A}$, where ${A\in \Q}$. 
\end{abstract}


\maketitle

\section{The Result}

Let ${\H=\{z\in \C:\Im z>0\}}$  be the Poincaré plane and~$j$ the $j$-invariant. 
The numbers of the form $j(\tau)$, where ${\tau\in \H}$ is an imaginary quadratic number, are called  \textsl{singular moduli}. It is known that $j(\tau)$ is an algebraic integer satisfying
$$
[\Q(\tau,j(\tau)):\Q(\tau)]= [\Q(j(\tau)):\Q]=h(\Delta),
$$
where~$\Delta$ is the discriminant of the complex multiplication  order ${\OO=\End\langle\tau,1\rangle}$ (the endomorphism ring of the lattice generated by~$\tau$ and~$1$) and ${h(\Delta)=h(\OO)}$ is the class number.

Let ${F(x_1,x_2)\in \C[x_1,x_2]}$ be an irreducible complex polynomial with 
$$
\deg_{x_1}F+\deg_{x_2}F\ge 2.
$$
In 1998 André~\cite{An98} proved that the equation ${F(j(\tau_1),j(\tau_2))=0}$ has at most finitely many solutions in singular moduli $j(\tau_1),j(\tau_2)$, unless $F(x_1,x_2)$ is the classical modular polynomial $\Phi_N(x_1,x_2)$ of some level~$N$. Recall that $\Phi_N$ is defined (up to a constant multiple) as the irreducible polynomial satisfying ${\Phi_N(j,j_N)=0}$, where ${j_N(z)=j(Nz)}$.

André's result was the first non-trivial contribution to the celebrated André-Oort conjecture on the special subvarieties of Shimura varieties; see~\cite{KY14,Pi11} and the references therein. 

Independently of André the same result was also obtained by Edixhoven~\cite{Ed98}, but Edixhoven had to assume the Generalized Riemann Hypothesis for certain $L$-series to be true. See also the work of Breuer~\cite{Br01}, who gave an explicit version of Edixhoven's result. 

Further proof followed; we mention specially the remarkable argument of Pila~\cite{Pi09}. It is based on an idea of Pila and Zannier~\cite{PZ08} and readily extends to higher dimensions~\cite{Pi11}.

The arguments of André and Pila are  non-effective, because they use the Siegel-Brauer lower bound for the class number.  
Recently Kühne~\cite{Ku12,Ku13} and, independently, Bilu, Masser, and Zannier~\cite{BMZ13} found unconditional effective proofs of André's theorem. Besides giving general results, both articles~\cite{Ku13} and~\cite{BMZ13} treat also some particular curves, showing they have no CM-points at all. For instance, Kühne \cite[Theorem~5]{Ku13}  proved  that a sum of two singular moduli can never be~$1$:
\begin{equation}
\label{esum1}
j(\tau_1)+j(\tau_2)\ne 1.
\end{equation}
Neither can their product be~$1$, as shown in~\cite{BMZ13}:
\begin{equation}
\label{eprod1}
j(\tau_1)j(\tau_2)\ne 1.
\end{equation}

A vast generalization of~\eqref{esum1} was given in~\cite{ABP14}: it is shown that, with ``obvious'' exceptions, two distinct singular moduli cannot satisfy a linear relation over~$\Q$. 

In this note we  combine ideas from~\cite{ABP14} and~\cite{BMZ13} and generalize~\eqref{eprod1}, showing that, with ``obvious'' exceptions, the product of two singular moduli cannot be a non-zero rational number. 

\begin{theorem}
\label{thprod}
Assume that ${j(\tau_1)j(\tau_2)\in \Q^\times}$. Then we have one of the following options:
\begin{itemize}
\item
(rational case)  both $j(\tau_1)$ and $j(\tau_2)$ are rational numbers (in fact integers); 
\item
(quadratic case)
$j(\tau_1)$ and $j(\tau_2)$ are of degree~$2$ and conjugate over~$\Q$. 
\end{itemize}
\end{theorem}

One may remark that the corresponding discriminants satisfy 
\begin{equation}
\label{edeg1}
h(\Delta_1)=h(\Delta_2)=1
\end{equation} 
in the rational  case and 
\begin{equation}
\Delta_1=\Delta_2=\Delta, \qquad h(\Delta)=2
\end{equation}
in the quadratic case. 

The full lists of discriminants of class numbers~$1$ and~$2$ are well-known. In particular, there are 13 discriminants of class number~$1$ (of which discriminant~$-3$ must be disregarded, because the corresponding $j$-value is~$0$) and 29 discriminants of class number~$2$. The former are reproduced in Table~\ref{talist1}, together with the corresponding $j$-invariants.
The latter are reproduced in Table~\ref{taquad},  together with the corresponding Hilbert Class Polynomials, calculated with \textsf{Sage}. (The contents of the right column of Table~\ref{taquad} is not relevant now; it will be used in Subsection~\ref{ssh2}.) 
This implies that there are 78 (unordered) pairs of the rational type and 29 pairs of the quadratic type.

\begin{table}
\caption{Discriminants~$\Delta$ with ${h(\Delta)=1}$ and the corresponding $j$-invariants}
\label{talist1}
{\footnotesize
\begin{equation*}
\begin{array}l
\begin{array}{r|lllllllll}
\Delta
&-3&-4&-7&-8&-11&-12&-16&-19&-27\\
j&0&1728&-3375&8000& -32768& 54000& 287496& -884736& -12288000\\
\end{array}\\
\hline
\begin{array}{r|lllllllll}
\Delta&-28&-43&-67&-163\\
j&16581375&-884736000&
-147197952000& -262537412640768000
\end{array}
\end{array}
\end{equation*}  
}
\end{table}

\begin{table}[t]
\caption{Discriminants~$\Delta$ with ${h(\Delta)=2}$  and their Hilbert Class Polynomials}
\label{taquad}
{\tiny
$$
\begin{array}{r|l|l}
\Delta& H_\Delta(x)=x^2+a_1x+a_0&a_0/a_1^2\\
\hline
-15& x^2 + 191025 x - 121287375 &-3.32\cdot10^{-3} \\
-20& x^2 - 1264000 x - 681472000 &-4.27\cdot10^{-4} \\
-24& x^2 - 4834944 x + 14670139392 &6.28\cdot10^{-4} \\
-32& x^2 - 52250000 x + 12167000000 &4.46\cdot10^{-6} \\
-35& x^2 + 117964800 x - 134217728000 &-9.65\cdot10^{-6} \\
-36& x^2 - 153542016 x - 1790957481984 &-7.60\cdot10^{-5} \\
-40& x^2 - 425692800 x + 9103145472000 &5.02\cdot10^{-5} \\
-48& x^2 - 2835810000 x + 6549518250000 &8.14\cdot10^{-7} \\
-51& x^2 + 5541101568 x + 6262062317568 &2.04\cdot10^{-7} \\
-52& x^2 - 6896880000 x - 567663552000000 &-1.19\cdot10^{-5} \\
-60& x^2 - 37018076625 x + 153173312762625 &1.12\cdot10^{-7} \\
-64& x^2 - 82226316240 x - 7367066619912 &-1.09\cdot10^{-9} \\
-72& x^2 - 377674768000 x + 232381513792000000 &1.63\cdot10^{-6} \\
-75& x^2 + 654403829760 x + 5209253090426880 &1.22\cdot10^{-8} \\
-88& x^2 - 6294842640000 x + 15798135578688000000 &3.99\cdot10^{-7} \\
-91& x^2 + 10359073013760 x - 3845689020776448 &-3.58\cdot10^{-11} \\
-99& x^2 + 37616060956672 x - 56171326053810176 &-3.97\cdot10^{-11} \\
-100& x^2 - 44031499226496 x - 292143758886942437376 &-1.51\cdot10^{-7} \\
-112& x^2 - 274917323970000 x + 1337635747140890625 &1.77\cdot10^{-11} \\
-115& x^2 + 427864611225600 x + 130231327260672000 &7.11\cdot10^{-13} \\
-123& x^2 + 1354146840576000 x + 148809594175488000000 &8.12\cdot10^{-11} \\
-147& x^2 + 34848505552896000 x + 11356800389480448000000 &9.35\cdot10^{-12} \\
-148& x^2 - 39660183801072000 x - 7898242515936467904000000 &-5.02\cdot10^{-9} \\
-187& x^2 + 4545336381788160000 x - 3845689020776448000000 &-1.86\cdot10^{-16} \\
-232& x^2 - 604729957849891344000 x + 14871070713157137145512000000000 &4.07\cdot10^{-11} \\
-235& x^2 + 823177419449425920000 x + 11946621170462723407872000 &1.76\cdot10^{-17} \\
-267& x^2 + 19683091854079488000000 x + 531429662672621376897024000000 &1.37\cdot10^{-15} \\
-403& x^2 + 2452811389229331391979520000 x - 108844203402491055833088000000 &-1.81\cdot10^{-26} \\
-427& x^2 + 15611455512523783919812608000 x + 155041756222618916546936832000000&6.36\cdot10^{-25}
\end{array}
$$}
\end{table}

\begin{remark}
Two pairs of non-zero rational singular moduli have the same product:
$$
1728\cdot (-147197952000)= 287496\cdot(-884736000)= -254358061056000. 
$$
All other products are pairwise distinct (and distinct from $-254358061056000$). 
Hence there exist $77$  products of two non-zero rational singular moduli. 
Denote by~$S$ the set consisting of these $77$ numbers,  and of
the free terms of the $29$ Hilbert class polynomials of the discriminants with ${h=2}$ (displayed in the central column of Table~\ref{taquad}). 
A quick verification shows that all these $106$ numbers are pairwise distinct, so the set~$S$ consists of exactly $106$ elements. Now Theorem~\ref{thprod} implies 
 that for every ${A\in S\smallsetminus \{-254358061056000\}}$ there exists exactly one unordered pair ${\{j(\tau_1),j(\tau_2)\}}$ of singular moduli such that ${j(\tau_1)j(\tau_2)=A}$, that for ${A=-254358061056000}$ there are exactly two such pairs, and 
for every non-zero rational ${A\notin S}$ there is no such pair at all.  Thus, we obtain a very explicit version of André's theorem for the one-parametric family  of hyperbolas ${x_1x_2=A}$, where ${A\in \Q^\times}$. 
\end{remark}

{\sloppy

\paragraph{Acknowledgments}
Yuri Bilu was supported by  the \textsl{Agence National de la Recherche} project ``Hamot'' (ANR 2010 BLAN-0115-01) and by the ALGANT scholarship program. While working on this project, he was enjoying the  hospitality of the Max-Plank-Institute of Mathematics, the University of Stellenbosch and the University of the Witwatersrand. 

}

We thank Bill Allombert, Alexander Boritchev, Florian Breuer and Lars Kühne for useful discussions, and the referee for helpful suggestions.

Our calculations were performed using the computer packages \textsf{PARI/GP}~\cite{PARI} and \textsf{Sage}~\cite{Sage}.

\section{Auxiliary Material}

\subsection{Estimates for the $j$-Invariant}

We denote by~$\gerD$ the standard fundamental domain: the open hyperbolic triangle with vertices 
$$
\zeta_3=\frac{-1+\sqrt{-3}}2, \quad \zeta_6=\frac{1+\sqrt{-3}}2, \quad \infty,
$$ 
together with the geodesics connecting~$\zeta_6$ with $\sqrt{-1}$ and with $\infty$. We write ${\gerD=\gerD^+\cup\gerD^-}$, where 
$$
\gerD^+=\{z\in \gerD: \Re z\ge 0\}, \qquad \gerD^-=\{z\in \gerD: \Re z< 0\}. 
$$
For ${z\in \H}$ we denote ${q_z=e^{2\pi \sqrt{-1}z}}$. 

When $j(z)$ is large, it is approximately of the same magnitude as $q_z^{-1}$. The following is Lemma~1 from~\cite{BMZ13}, which makes this explicit. 

\begin{proposition}
\label{pr2079}
For ${z\in \gerD}$ we have ${\bigl||j(z)|-|q_z^{-1}|\bigr|\le2079}$.
\end{proposition}

Next let us study how small  can $j(z)$ be. Clearly, if  ${z\in \gerD^+}$ is such that $|j(z)|$ is small then~$z$ must be close to~$\zeta_6$ and $|j(z)|$ must be of magnitude ${|z-\zeta_6|^3}$, because~$j$ has a triple zero  at~$\zeta_6$. We again want to make this explicit. 

\begin{proposition}
\label{prjsmall}
For  ${z\in \gerD^+}$ one of the following alternatives holds:
when ${|z-\zeta_6|\ge10^{-3}}$ we have  ${|j(z)|\ge 4.4\cdot10^{-5}}$, and when  ${|z-\zeta_6|\le10^{-3}}$ we have
\begin{equation}
\label{enear}
44000|z-\zeta_6|^3\le |j(z)|\le 47000|z-\zeta_6|^3. 
\end{equation}
\end{proposition}

\begin{remark}
\begin{enumerate}
\item
The same statement (with the same proof) holds true for ${z\in \gerD^-}$ with~$\zeta_6$ replaced by~$\zeta_3$.

\item
Only the lower bound from~\eqref{enear} will be used in the sequel.

\end{enumerate}
\end{remark}

The proof of Proposition~\ref{prjsmall} requires a   Schwarz-type lemma.

\begin{lemma}
\label{lschw}
Let~$f$ be a holomorphic function in an open neighborhood of the disc ${|z-a|\le R}$ and assume that ${|f(z)|\le B}$ in this disc. Further, let~$\ell$ be a non-negative integer such that ${f^{(k)}(a)=0}$ for ${0\le k<\ell}$ and ${f^{(\ell)}(a)\ne 0}$. 
Set ${A=f^{(\ell)}(a)/\ell!}$\,. Then in the same disc  we have the estimate
\begin{equation*}
\left|f(z)-A(z-a)^\ell\right|\le \frac{|A|R^\ell+B}{R^{\ell+1}}|z-a|^{\ell+1}. 
\end{equation*}
\end{lemma}

\paragraph{Proof}
The function ${g(z)=(f(z)-A(z-a)^\ell)(z-a)^{-\ell-1}}$ is holomorphic in an open neighborhood of the disc ${|z-a|\le R}$, and on the circle ${|z-a|=R}$ we have the estimate
\begin{equation*}
|g(z)|\le \frac{|A|R^\ell+B}{R^{\ell+1}}
\end{equation*}
By the maximal principle the same estimate holds true in the disc ${|z-a|\le R}$. Whence the result. \qed

\paragraph{Proof of Proposition~\ref{prjsmall}}
Our starting point is the estimate ${|j(z)|\le 23000}$ in the disc ${|z-\zeta_6|\le {\sqrt3}/4}$, see Lemma~2 in~\cite{BMZ13}. (The statement of the lemma has~$30000$, but the actually proved upper bound is~$23000$.) We also have ${j(\zeta_6)=j'(\zeta_6)=j''(\zeta_6)=0}$ and 
$$
j'''(\zeta_6)=-162\Gamma(1/3)^{18}\pi^{-9}\sqrt{-1}=-\sqrt{-1}\cdot274470.48\dots
$$
(see, for instance, \cite[page 777]{Wu14}). Using Lemma~\ref{lschw} with 
$$
A=j'''(\zeta_6)/6=-\sqrt{-1}\cdot45745.08\dots, \quad R= {\sqrt3}/4, \quad B=23000,
$$
we obtain in the disc ${|z-\zeta_6|\le {\sqrt3}/4}$ the estimate
\begin{equation}
\label{eschwj}
\left|j(z)-A(z-\zeta_6)^3\right|\le \frac{46000({\sqrt3}/4)^3+23000}{({\sqrt3}/4)^4}|z-\zeta_6|^4<761000|z-\zeta_6|^4. 
\end{equation}
In the disc ${|z-\zeta_6|\le10^{-3}}$ the right-hand side of~\eqref{eschwj} is bounded by ${761|z-\zeta_6|^3}$, and we obtain
$$
(|A|-761)|z-\zeta_6|^3\le |j(z)|\le (|A|+761)|z-\zeta_6|^3,
$$
which proves~\eqref{enear}.  

In particular, on the circle ${|z-\zeta_6|= 10^{-3}}$ we have  ${|j(z)|\ge 4.4\cdot10^{-5}}$. 
Using the properties of~$j$ on the boundary of~$\gerD^+$ (where it takes real values), we deduce from this that the estimate ${|j(z)|\ge 4.4\cdot10^{-5}}$ holds for any~$z$ on the boundary of the domain
${\gerD^+\cap \{z: |z-\zeta_6|\ge 10^{-3}\}}$. By the maximum principle this is true for any~$z$ in this domain, which completes the proof. \qed

\subsection{The Conjugates of $j(\tau)$}
\label{ssconj}
Let ${\tau\in \H}$ be imaginary quadratic,  and~$\Delta$ be the discriminant of its CM-order.    It is well-known that the $\Q$-conjugates of the algebraic integer $j(\tau)$ can be described explicitly. Below we briefly recall this description.

Denote by ${T=T_\Delta}$ the set of triples of integers $(a,b,c)$ such that 
\begin{equation*}
\begin{gathered}
\gcd(a,b,c)=1, \quad \Delta=b^2-4ac,\\
\text{either\quad $-a < b \le a < c$\quad or\quad $0 \le b \le a = c$}
\end{gathered}
\end{equation*}
\begin{proposition}
\label{pallconj}
The map
\begin{equation}
\label{ejabcmap}
(a,b,c)\mapsto j\left(\frac{b+\sqrt{\Delta}}{2a}\right)
\end{equation}
defines a bijection from $T_\Delta$ onto the set of $\Q$-conjugates
of $j(\tau)$. 
In particular, ${h(\Delta)=|T_\Delta|}$.
\end{proposition}

In the sequel it will be convenient to use the notation 
$$
\tau(a,b,c)=\frac{b+\sqrt{\Delta}}{2a}
$$
for ${(a,b,c)\in T_\Delta}$. One may notice that $\tau(a,b,c)$   belongs to the standard fundamental domain~$\gerD$.

\paragraph{Proof}
Let~$\OO$ be the imaginary quadratic order of discriminant~$\Delta$ and $\Cl(\OO)$ its class group. Then the set of conjugates of $j(\tau)$ coincides with the set 
${\{j(\calA): \calA\in \Cl(\OO)\}}$, see \cite[Proposition 13.2]{Co89}.

On the other hand, for every ${(a,b,c)\in T_\Delta}$ the lattice ${\langle1,\tau(a,b,c)\rangle}$ is an invertible fractional ideal of~$\OO$. Denote by $\calA(a,b,c)$ the class of this ideal. Clearly, ${j(\tau(a,b,c))=j(\calA(a,b,c))}$. 
Finally,  joint application of Theorems~2.8 and~7.7 from~\cite{Co89} implies that the map 
${(a,b,c)\mapsto \calA(a,b,c)}$
defines a bijection of the sets $T_\Delta$ and $\Cl(\OO)$. This completes the proof of Proposition~\ref{pallconj}. \qed

\bigskip

Another useful observation: since ${0<a\le c}$ and ${|b|\le a}$, we have 
\begin{equation}
\label{eadelta}
|\Delta|=4ac-b^2\ge 4a^2-a^2=3a^2. 
\end{equation}

The following statement is proved by a straightforward verification; we omit the details. 

\begin{proposition}
\label{pra12}
For every negative discriminant~$\Delta$ the set~$T_\Delta$ has exactly one element $(a,b,c)$ with ${a=1}$ and at most two elements with ${a=2}$. More precisely,~$T_\Delta$ has:
\begin{itemize}
\item
two elements with ${a=2}$ if ${\Delta\equiv 1\bmod 8}$, ${\Delta\ne -7}$;

\item
one element with ${a=2}$ if ${\Delta\equiv 8,12\bmod 16}$, ${\Delta\ne -4,-8}$;

\item
no elements with ${a=2}$ if ${\Delta\equiv 5\bmod 8}$ or ${\Delta\equiv 0,4\bmod 16}$. \qed

\end{itemize}

\end{proposition}

\subsection{Hilbert Class Polynomials}
\label{sshilb}
The monic polynomial having the numbers on the right of~\eqref{ejabcmap} as roots is usually called the \textsl{Hilbert class polynomial} of discriminant~$\Delta$; we denote it by $H_\Delta(x)$. It is a  polynomial in $\Z[x]$ of degree ${h=h(\Delta)}$.

If ${j(\tau)j(\tau')=A\in \Q^\times}$, and ${\Delta,\Delta'}$ are the corresponding discriminants, then ${h(\Delta)=h(\Delta')=h}$; furthermore, writing 
$$
H_\Delta(x)=x^h+a_{h-1}x^{h-1}+\cdots+a_0, \quad H_{\Delta'}(x)=x^h+a_{h-1}'x^{h-1}+\cdots+a_0',
$$ 
we have ${a'_{h-i} =A^ia_i/a_0}$ for ${i=0, \ldots,h}$, where we set ${a_h=a_h'=1}$. It follows that
\begin{equation}
\label{ehilbprod}
\frac{a_{i-1}a_{i+1}}{a_i^2}=\frac{a_{h-i-1}'a_{h-i+1}'}{(a_{h-i}')^2} \qquad (i=1, \ldots, h-1)
\end{equation}
with the obvious convention in the case of zero denominators (that is, if the denominator on the left is zero, but the numerator isn't, then we have the same on the right). In particular, if ${\Delta=\Delta'}$ then 
\begin{equation}
\label{ehilbprodequ}
\frac{a_{i-1}a_{i+1}}{a_i^2}=\frac{a_{h-i-1}a_{h-i+1}}{a_{h-i}^2} \qquad (i=1, \ldots, h-1).
\end{equation}
This gives
an easy way to exclude the possibility that 
${j(\tau)j(\tau')\in \Q^\times}$ in every concrete case. 

\subsection{Comparing Two CM Fields}

The following theorem is proved in~\cite{ABP14}, see Corollary~4.2 and Proposition~4.3.

{\sloppy

\begin{theorem}
\label{thedan}
Let~$\tau_1$ and~$\tau_2$ be imaginary quadratic numbers such that  ${\Q(j(\tau_1))=\Q(j(\tau_2))}$, and let ${\Delta_1,\Delta_2}$  
be the  discriminants of the corresponding CM orders.  

\begin{enumerate}

\item
\label{ine}
Assume that  ${\Q(\tau_1)\ne \Q(\tau_2)}$.
Then the field ${L=\Q(j(\tau_1))=\Q(j(\tau_2))}$ is one of the fields in Table~\ref{tadouble}. 
\begin{table}[t]
\caption{Fields presented as $\Q(j(\tau_1))$ and $\Q(j(\tau_2))$ with ${\Q(\tau_1)\ne \Q(\tau_2)}$}
\label{tadouble}
{\footnotesize
$$
\begin{array}{l|l|l}
\text{Field $L$}&h=[L:\Q]&\text{discriminants~$\Delta$ of CM-orders $\End\langle1,\tau\rangle$}\\
&& \text{such that ${L=\Q(j(\tau))}$}
\\
\hline
\Q&1& -3, -4, -7, -8, -11, -12, -16, -19, -27, -28, -43, \\
&&-67, -163
\\
\hline
\Q(\sqrt{2})&2&-24, -32, -64, -88
\\
\Q(\sqrt{3})&2&-36, -48
\\
\Q(\sqrt{5})&2&-15, -20, -35, -40, -60, -75, -100, -115, -235
\\
\Q(\sqrt{13})&2&-52, -91, -403
\\
\Q(\sqrt{17})&2&-51, -187
\\
\hline
\Q(\sqrt{2},\sqrt{3})&4&-96, -192, -288
\\
\Q(\sqrt{3},\sqrt{5})&4&-180, -240
\\
\Q(\sqrt{5},\sqrt{13})&4&-195, -520, -715
\\
\Q(\sqrt{2},\sqrt{5})&4&-120, -160, -280, -760
\\
\Q(\sqrt{5},\sqrt{17})&4&-340, -595
\\
\hline
\Q(\sqrt{2},\sqrt{3},\sqrt{5})&8&-480, -960
\end{array}\\
$$
}\end{table}

\item
\label{ieq}
Assume that ${\Q(\tau_1)= \Q(\tau_2)}$. Then either ${\Delta_1, \Delta_2 \in \{-3,-12,-27\}}$ or ${\Delta_1/\Delta_2\in \{1,4,1/4\}}$. 
\end{enumerate}
\end{theorem}

}

\section{Proof of Theorem~\ref{thprod}}

Assume that ${j(\tau_1)j(\tau_2)=A\in \Q^\times}$. In fact, since both $j(\tau_i)$ are algebraic integers, we must have ${A\in \Z}$, but this will not play any significant role in the sequel. As before, we denote by~$\Delta_1$ and~$\Delta_2$ the discriminants of the corresponding CM orders. 

First of all, we clearly have ${\Q(j(\tau_1))=\Q(j(\tau_2))}$, and, in particular, 
$$
h(\Delta_1)=h(\Delta_2)=h. 
$$
When ${h=1}$ we are, obviously, in the ``rational case'' of Theorem~\ref{thprod}, and 
when ${h=2}$, we are in the ``quadratic case''; this is slightly less obvious, see Subsection~\ref{ssh2}.

Therefore in the sequel we will assume that ${h\ge3}$, and, in particular,
$$
|\Delta_1|,|\Delta_2|\ge 23. 
$$
We also have ${j(\tau_1)\ne j(\tau_2)}$, because otherwise ${j(\tau_1)^2\in \Q}$, which implies ${h\le 2}$.

We will bound~$A$ from below and from above in terms of~$\Delta_1$ and~$\Delta_2$,  and will see that  the two bounds contradict each other in all but finitely many cases. These remaining few cases can be treated by direct verification using any available number-theoretic computer package (we used \textsf{PARI}~\cite{PARI}). 

\subsection{Lower Bound}
We want to bound~$A$ from below. We may assume that, for ${i=1,2}$, there exist triples ${(a_i,b_i,c_i)\in T_{\Delta_i}}$  such that ${\tau_i=\tau(a_i,b_i,c_i)}$; see Subsection~\ref{ssconj} for the details.    Proposition~\ref{pra12} implies that,  conjugating over~$\Q$, we may assume that ${a_1=1}$. Then ${q_{\tau_1}=e^{2\pi \sqrt{-1}\tau_1}=\pm e^{-\pi|\Delta_1|^{1/2}}}$. Using Proposition~\ref{pr2079}, we obtain
\begin{equation}
\label{ejt1}
|j(\tau_1)|\ge e^{\pi|\Delta_1|^{1/2}}-2079 \ge 0.9994e^{\pi|\Delta_1|^{1/2}},
\end{equation}
the latter estimate being valid because  ${|\Delta_1|\ge 23}$. 

To bound $|j(\tau_2)|$ from below we use Proposition~\ref{prjsmall}. We may assume that~$\tau_2$ belongs to the right half~$\gerD^+$ of the fundamental domain, the case ${\tau_2\in\gerD^-}$ being absolutely analogous\footnote{In fact~$\tau_2$ does belong to $\gerD^+$, and even to the boundary of~$\gerD^+$ (because $j(\tau_1)$ is real and so is $j(\tau_2)$) but this is of no importance for us.}. 

Proposition~\ref{prjsmall} implies that 
$$
|j(\tau_2)|\ge \min\{4.4\cdot10^{-5}, 44000|\tau_2-\zeta|^3\}, 
$$
where ${\zeta=\zeta_6}$. 
To estimate ${|\tau_2-\zeta|}$ we first notice that, since ${j(\tau_2)\ne 0}$ we have ${\tau_2\ne \zeta}$, and since ${\tau_2\in \gerD^+}$, we have
$$
\frac{\sqrt{|\Delta_2|}}{2a_2}=\Im\tau_2>\frac{\sqrt3}{2}=\Im\zeta. 
$$
Therefore
{\small
$$
|\tau_2-\zeta|\ge \left|\frac{\sqrt{|\Delta_2|}}{2a_2}-\frac{\sqrt3}{2}\right|= \frac{\bigl||\Delta_2|-3a_2^2\bigr|}{2a_2(\sqrt{|\Delta_2|}+a_2\sqrt3)}\ge \frac{1}{2a_2(\sqrt{|\Delta_2|}+a_2\sqrt3)}\ge \frac{\sqrt3}{4|\Delta_2|}, 
$$}%
the last inequality being implied by~\eqref{eadelta}. We obtain 
$$
|j(\tau_2)|\ge \min\{4.4\cdot10^{-5}, 3500|\Delta_2|^{-3}\}.
$$
Combined with~\eqref{ejt1}, this results in the following lower estimate:
\begin{equation}
\label{elowereq}
|A|\ge 3000e^{\pi|\Delta_1|^{1/2}}\min\{10^{-8}, |\Delta_2|^{-3}\}. 
\end{equation}

\subsection{Upper Bound}
Now we want to bound~$A$ from above. It will be convenient to consider separately the cases ${\Delta_1=\Delta_2}$ and ${\Delta_1\ne \Delta_2}$, and for the latter also separate the sub-cases ${\Q(\tau_1)=\Q(\tau_2)}$ and ${\Q(\tau_1)\ne\Q(\tau_2)}$. The arguments differ only technically, therefore we give the full details only in the first of the three cases.

\subsubsection{The Case ${\Delta_1=\Delta_2=\Delta}$}


In this case the singular moduli $j(\tau_1)$ and $j(\tau_2)$ are conjugate over~$\Q$. Since ${j(\tau_1)\ne j(\tau_2)}$ and ${j(\tau_1)j(\tau_2)\in \Q}$, the field ${L=\Q(j(\tau_1))=\Q(j(\tau_1))}$ admits a non-trivial automorphism of order~$2$, swapping $j(\tau_1)$ and $j(\tau_2)$. It follows that ${h=[L:\Q]}$ is an even number; in particular, ${h\ge 4}$.

To bound~$A$ from above, we must impose certain assumptions on our  discriminants. We obtain two bounds: one sharper, under a more restrictive assumption, the other less
sharp, but valid under a milder assumption. Precisely:

\begin{enumerate}
\item 
\label{i23}
assume that ${|\Delta|\ge103}$, 
that ${h>4}$ when ${\Delta\equiv 8,12\bmod 16}$ and that ${h>6}$ when ${\Delta\equiv 1 \bmod8}$; 
then we have ${|A|\le 1.11e^{(2\pi/3)|\Delta|^{1/2}}}$;

\item
\label{i56}
assume  that ${|\Delta|\ge399}$ and 
that ${h>4 }$  when ${\Delta\equiv 1 \bmod8}$; 
then we have ${|A|\le 1.001e^{(5\pi/6)|\Delta|^{1/2}}}$.

\end{enumerate}

Since both $j(\tau_1)$ and $j(\tau_2)$ generate the same field of degree~$h$ over~$\Q$, the Galois orbit of the pair $(j(\tau_1), j(\tau_2))$ (over~$\Q$) has exactly~$h$ elements; moreover, each conjugate of $j(\tau_1)$ occurs exactly once as the first coordinate of a pair in the orbit, and each conjugate of $j(\tau_2)$ occurs exactly once as the second coordinate. Every such conjugate pair is of form ${(j(\tau_1'),j(\tau_2'))}$, where, for ${i=1,2}$ we have ${\tau_i'=\tau(a_i,b_i,c_i)}$ for some triples ${(a_i,b_i,c_i)\in T_\Delta}$ (see Subsection~\ref{ssconj}). 

For the proof of item~\ref{i23}  call the pair ${(j(\tau_1'),j(\tau_2'))}$ ``good'' if ${a_1,a_2\ge3}$ and ``bad'' otherwise. Proposition~\ref{pra12} implies that there are at most~$6$ ``bad'' pairs in the case ${\Delta\equiv 1 \bmod8}$, at most~$4$ ``bad'' pairs in the case ${\Delta\equiv 8,12\bmod 16}$,  and at most~$2$ ``bad'' pairs in all other cases. Hence at least one ``good'' pair exists in any case (recall that ${h\ge 4}$).   For such a pair we have
 $$
|q_{\tau_i}|= e^{(\pi/a_i)|\Delta|^{1/2}}\le e^{(\pi/3)|\Delta|^{1/2}}. 
$$
Now using Proposition~\ref{pr2079}, we obtain
$$
|A|=|j(\tau_1')j(\tau_2')|\le \bigl(e^{(\pi/3)|\Delta|^{1/2}}+2079\bigr)^2,
$$
which is bounded above by ${1.11e^{(2\pi/3)|\Delta|^{1/2}}}$ because ${|\Delta|\ge 103}$. 

For the proof of item~\ref{i56}  call a pair ``good'' if  ${a_1,a_2\ge2}$ and ${a_1+a_2\ge 5}$. The rest of the proof is  similar to that of item~\ref{i23}, and we omit the details. 


\bigskip

Now, combining~\eqref{elowereq} with the bound ${|A|\le 1.11e^{(2\pi/3)|\Delta|^{1/2}}}$ yields ${|\Delta|<103}$, and combining it with ${|A|\le 1.01e^{(5\pi/6)|\Delta|^{1/2}}}$  yields ${|\Delta|<399}$. This shows that~$\Delta$ must satisfy  one of the following conditions:

\begin{enumerate}
\item
\label{icond1}
\quad ${h(\Delta)\ge 4}$ is even and ${|\Delta|<103}$;

\item
\label{icond2}
\quad
${\Delta\equiv 8,12\bmod 16}$, \ ${h(\Delta)=4}$ and ${103\le|\Delta|<399}$; 


\item
\label{icond3}
\quad 
${\Delta\equiv 1 \bmod8}$, \ ${h(\Delta)=6}$ and ${103\le|\Delta|<399}$; 

\item
\label{icond4}
\quad ${\Delta\equiv 1 \bmod8}$, \  ${h(\Delta)=4}$ and  ${|\Delta|\ge 103}$.

\end{enumerate}

There are not too many~$\Delta$ satisfying one of these: no discriminant satisfies condition~\ref{icond4}, and the full lists (found with \textsf{PARI}) for conditions~\ref{icond1},~\ref{icond2}  and~\ref{icond3} are, respectively,
\begin{align*}
&-63, -80, -96, -39, -55, -56, -68, -84, -87, -95;\\
&-196, -180, -132, -228, -292, -340, -372, -388;\\
&-175, -135, -207, -247.
\end{align*}
Using \textsf{PARI}, we computed the Hilbert Class Polynomials for the discriminants above, and verified that none of these polynomials satisfies condition~\eqref{ehilbprodequ}. We do not include the results of this computation in the article, but they can be obtained from the authors, together with the source codes.

\subsubsection{The Case ${\Delta_1\ne \Delta_2}$, but ${\Q(\tau_1)=\Q(\tau_2)}$}
Theorem~\ref{thedan} implies that we have one of the options ${\Delta_1, \Delta_2 \in \{-3,-12,-27\}}$ or ${\Delta_1/\Delta_2\in \{4,1/4\}}$. In the former case  ${h=1}$, which is excluded. Hence, we may assume that ${\Delta_1=4\Delta_2}$ and write
${\Delta_1=4\Delta}$, ${\Delta_2=\Delta}$.

Observe that ${\Delta\equiv 1 \bmod 8}$. Indeed, recall the ``class number formula'' 
$$
h(m^2\Delta)= \frac{m}{\omega}\prod_{p\mid m}\left(1-\frac1p\left(\frac{\Delta}{p}\right)\right) h(\Delta), \qquad 
\omega=
\begin{cases}
3,&\Delta=-3,\\
2,&\Delta=-4,\\
1,&\text{otherwise},
\end{cases}
$$
where $(\Delta/p)$ is the Kronecker symbol (see \cite[Corollary~7.28]{Co89}). In our case ${\omega=1}$ and ${m=2}$, which gives 
${h(4\Delta)=(2-(\Delta/2))h(\Delta)}$, and for the equality ${h(\Delta)=h(4\Delta)}$ we must have ${(\Delta/2)=1}$, that is, ${\Delta\equiv 1 \bmod 8}$. 


The lower bound~\eqref{elowereq} becomes
$$
|A|\ge 3000e^{2\pi|\Delta|^{1/2}}\min\{10^{-8}, |\Delta|^{-3}\}.
$$
For the upper bound we argue as in the previous subsection. We use a pair ${(j(\tau_1'),j(\tau_2'))}$ with ${a_1,a_2\ge 2}$; this is always possible because ${h\ge 3}$. Since 
$$
\Delta_1=4\Delta\equiv4\bmod 32, 
$$
Proposition~\ref{pra12} implies that ${a_1\ge 3}$, which gives the upper bound 
$$
|A|\le (e^{(2\pi/3)|\Delta|^{1/2}}+2079)(e^{(\pi/2)|\Delta|^{1/2}}+2079) \le 2.4e^{(7\pi/6)|\Delta|^{1/2}}
$$ 
(we again use ${|\Delta|\ge 23}$). Comparing the two bounds, we deduce ${|\Delta|<20}$, contradicting the assumption ${h(\Delta)\ge 3}$. This completes the proof in this case.

\subsubsection{The Case ${\Delta_1\ne \Delta_2}$ and  ${\Q(\tau_1)\ne\Q(\tau_2)}$}
\label{sssnene}
According to Theorem~\ref{thedan}, we are now in one of the cases featured in sections ${h=4}$ and ${h=8}$ of  Table~\ref{tadouble}. We may assume ${|\Delta_1|>|\Delta_2|}$, and, inspecting the table, we find ${|\Delta_2|\ge 96}$ and ${|\Delta_1|\ge 160}$.   To bound~$|A|$ from above, we proceed as in the previous subsections. We use a pair ${(j(\tau_1'),j(\tau_2'))}$ with ${a_1\ge 3}$ and ${a_2\ge 2}$; this is always possible because ${h\ge 4}$ and none of our discriminants is ${1\bmod8}$. We obtain the upper bound 
\begin{equation}
\label{euppne}
|A|\le (e^{(\pi/3)|\Delta_1|^{1/2}}+2079)(e^{(\pi/2)|\Delta_2|^{1/2}}+2079) \le 1.005e^{(\pi/3)|\Delta_1|^{1/2}+(\pi/2)|\Delta_2|^{1/2}}. 
\end{equation}
Comparing it with~\eqref{elowereq}, we obtain
\begin{equation}
\label{enene}
\frac{2\pi}3|\Delta_1|^{1/2} +\log\frac{3000}{1.005}\le \frac\pi2|\Delta_2|^{1/2}+\max\{8\log10,3\log|\Delta_2|\}. 
\end{equation}
\begin{table}[t]
\caption{Data for Subsection~\ref{sssnene}}
\label{tanene}
{\footnotesize
$$
\begin{array}{rll}
\Delta & \frac{2\pi}3|\Delta|^{1/2} +\log\frac{3000}{1.005} &
\frac\pi2|\Delta|^{1/2}+\max\{8\log10,3\log|\Delta|\}\\[1mm]
\hline
-96	&	28.52217731&	33.81127871\\
-192	&	37.02216985&	40.18627311\\
-288	&	43.54444353&	45.07797837\\
\hline			
-180	&	36.10063895&	39.49512494\\
-240	&	40.44760943&	42.7553528\\
\hline			
-195	&	37.24801598&	40.35565771\\
-520	&	55.76093655&	54.58115383\\
-715	&	64.00442418&    61.71913074\\
\hline		
-120	&	30.94931641&	35.62789237\\
-160	&	34.49860294&	38.28985728\\
-280	&	43.05230081&	44.70513067\\
-760	&	65.74485596&	63.2038216\\
\hline			
-340	&	46.62510508&	47.38473388\\
-595	&	59.09415527&	57.481525\\
\hline			
-480	&	53.89226524&	52.93578157\\
-960	&	72.89882638&	69.27014397
\end{array}
$$}
\end{table}%
As Table~\ref{tanene} shows, the only case when~\eqref{enene} is satisfied is when ${\Delta_1=-160}$ and ${\Delta_2=-120}$. However, in this case ${\Delta_1\equiv 0\mod 16}$, which implies that we can use the pair ${(j(\tau_1'),j(\tau_2'))}$ with ${a_1,a_2\ge 3}$. This allows one to replace~\eqref{euppne} by a sharper bound 
$$
|A|\le (e^{(\pi/3)|\Delta_1|^{1/2}}+2079)(e^{(\pi/3)|\Delta_2|^{1/2}}+2079).
$$
A quick calculation shows that this bound contradicts~\eqref{elowereq} when  ${\Delta_1=-160}$, ${\Delta_2=-120}$.

\subsection{The case ${h=2}$}
\label{ssh2}
In this case we have the following three options:
\begin{itemize}
\item
${j(\tau_1)=j(\tau_2)}$;
\item
${j(\tau_1)\ne j(\tau_2)}$, ${\Delta_1=\Delta_2}$;
\item
${\Delta_1\ne\Delta_2}$. 
\end{itemize}

The second option is exactly the ``quadratic case'' of Theorem~\ref{thprod}. We are left with showing that the other two options are impossible. We use the data from Table~\ref{taquad}. 

If ${j(\tau_1)=j(\tau_2)}$ then ${j(\tau_1)^2=A}$, which means that the Hilbert class polynomials $H_\Delta(x)$ for ${\Delta=\Delta_1=\Delta_2}$ must be ${x^2-A}$. However, all polynomials in the second column of Table~\ref{taquad} have their middle coefficient $a_1$ distinct from~$0$.

Finally, if  ${\Delta_1\ne\Delta_2}$, then the quantity appearing in the third column of Table~\ref{taquad} must be the same for these two discriminants, see identity~\eqref{ehilbprod}. However, all entries in this column are distinct. \qed

{\footnotesize

\end{document}